\documentclass[11pt]{article}

\usepackage{mathtools}
\usepackage{color}
\usepackage{microtype}
\usepackage[usenames,dvipsnames,svgnames,table]{xcolor}
\usepackage{amsfonts,amsmath, amssymb,amsthm,amscd,mathrsfs}
\usepackage[height=9in,width=6.5in]{geometry}
\usepackage{verbatim}
\usepackage{tikz}
\usetikzlibrary{decorations.markings}
\tikzset{->-/.style={decoration={
			markings,
			mark=at position .6 with {\arrow{>}}},postaction={decorate}}}
\usepackage{tkz-euclide}

\usepackage[margin=10pt,font=small,labelfont=bf]{caption}
\usepackage{latexsym}
\usepackage[pdfauthor={ },bookmarksnumbered,hyperfigures,colorlinks=true,citecolor=BrickRed,linkcolor=BrickRed,urlcolor=BrickRed,pdfstartview=FitH]{hyperref}

\usepackage{graphicx}
\usepackage{enumitem}

\usepackage[inline,nolabel]{showlabels}

\usepackage[capitalize]{cleveref}
\Crefname{fact}{Fact}{Facts}
\crefname{theorem}{Theorem}{Theorems}
\crefname{thm}{Theorem}{Theorems}
\crefname{lemma}{Lemma}{Lemmas}
\crefname{claim}{Claim}{Claims}
\crefname{lem}{Lemma}{Lemmas}
\crefname{remark}{Remark}{Remarks}
\crefname{prop}{Proposition}{Propositions}
\crefname{defn}{Definition}{Definitions}
\crefname{corollary}{Corollary}{Corollaries}
\crefname{conjecture}{Conjecture}{Conjectures}
\crefname{question}{Question}{Questions}
\crefname{chapter}{Chapter}{Chapters}
\crefname{section}{Section}{Sections}
\crefname{part}{Part}{Parts}
\crefname{figure}{Figure}{Figures}

\newtheorem{theorem}{Theorem}[section]
\newtheorem{corollary}[theorem]{Corollary}
\newtheorem{lemma}[theorem]{Lemma}

\newtheorem{example}[theorem]{Example}

\newtheorem{remark}[theorem]{Remark}

\numberwithin{equation}{section}

\newtheorem*{theorem*}{Theorem}
\newtheorem*{claim*}{Claim}
\newtheorem*{conj*}{Conjecture}
\newtheorem*{corollary*}{Corollary}
\newtheorem*{definition*}{Definition}
\newtheorem*{example*}{Example}
\newtheorem*{exercise*}{Exercise}
\newtheorem*{lemma*}{Lemma}
\newtheorem*{observation*}{Observation}
\newtheorem*{proposition*}{Proposition}
\newtheorem*{question*}{Question}
\newtheorem*{remark*}{Remark}

\DeclareMathSymbol{\leqslant}{\mathalpha}{AMSa}{"36} 
\DeclareMathSymbol{\geqslant}{\mathalpha}{AMSa}{"3E} 
\DeclareMathSymbol{\eset}{\mathalpha}{AMSb}{"3F}     
\renewcommand{\le}{\;\leqslant\;}                   
\renewcommand{\ge}{\;\geqslant\;}                   
\renewcommand{\leq}{\;\leqslant\;}                   

\newcommand{\be}{\begin{equation}}
	\newcommand{\ee}{\end{equation}}
\newcommand\ba{\begin{align}}
	\newcommand\ea{\end{align}}

\usepackage{nicematrix}
\usepackage{xcolor}  

\newcommand{\expect}{\mathbb E}
\newcommand{\prob}{\mathbb P}
\newcommand{\bx}{\boldsymbol x}
\newcommand{\bz}{\boldsymbol z}
\newcommand{\ind}{\mathbf 1}
\DeclareMathOperator{\Cov}{Cov}
\DeclareMathOperator{\rank}{rank}

\begin{document}
	
\author{
	Pengfei Tang\thanks{Center for Applied Mathematics and KL-AAGDM, Tianjin University, Tianjin, 300072, China.
		Email: \textsf{pengfei\_tang@tju.edu.cn}. Supported by the National Natural Science Foundation of China No. 12571151.}
	\qquad
	Zibo Zhang\thanks{School of Mathematics, Tianjin University, Tianjin, 300350, China.
		Email: \textsf{19932788533@163.com}.}
}
\date{\today}
\title{Pairwise correlation bounds for the random-cluster model with cluster weight at most one}

\maketitle


\begin{abstract}
We establish quantitative bounds on pairwise edge correlations for the random-cluster model with cluster weight $q\in[0,1]$ on arbitrary finite connected graphs. We give two explicit bounds that retain the dependence on the activities of the two specified edges and hold uniformly over the activities of all remaining edges. We also obtain a refinement of one of these bounds using connectivity information from the graph with the two specified edges removed. These results sharpen the bound obtained by a direct application of Lorentzian polynomial theory. At $q=0$, the model is supported on connected spanning subgraphs. In particular, for the uniform probability measure on these subgraphs, the joint inclusion probability of two distinct edges is at most $10/9$ times the product of their marginal inclusion probabilities.
\end{abstract}

\section{Introduction}
\label{sec:introduction}

Let $G=(V,E)$ be a finite connected undirected graph, with loops and parallel
edges allowed, and suppose that $m:=|E|\ge2$. Parallel edges are distinct
elements of $E$. For $A\subseteq E$, let $k(A)$ be the number of connected components
of $(V,A)$, including isolated vertices. Given activities
$\bx=(x_g)_{g\in E}\in(0,\infty)^E$ and $q\in[0,1]$, define
\begin{equation}
Z_G(q,\bx)=\sum_{A\subseteq E}q^{k(A)-1}\bx^A,
\qquad
\mu_{G,q,\bx}(A)=\frac{q^{k(A)-1}\bx^A}{Z_G(q,\bx)},
\label{eq:measure}
\end{equation}
where $\bx^A:=\prod_{g\in A}x_g$, the empty product is one, and $0^0:=1$.
For $q>0$, this is exactly the usual random-cluster measure. Indeed, if
\[
Z_G^{\mathrm{cl}}(q,\bx):=\sum_{A\subseteq E}q^{k(A)}\bx^A,
\]
then $Z_G^{\mathrm{cl}}(q,\bx)=qZ_G(q,\bx)$ and
$q^{k(A)}\bx^A/Z_G^{\mathrm{cl}}(q,\bx)
=\mu_{G,q,\bx}(A)$. We use the reduced partition function
$Z_G=q^{-1}Z_G^{\mathrm{cl}}$ so that the common vanishing factor $q$
is removed before taking $q\downarrow0$. This normalization changes no
probability or correlation for $q>0$ and gives a direct definition at
$q=0$, where
\[
Z_G(0,\bx)=\sum_{A:\,k(A)=1}\bx^A,\qquad
\mu_{G,0,\bx}(A)
=\frac{\ind_{\{k(A)=1\}}\bx^A}
{\displaystyle\sum_{B:\,k(B)=1}\bx^B}.
\]
Thus the endpoint measure is supported on connected spanning subgraphs,
with weights $\bx^A$. Since the sum in \eqref{eq:measure} is finite and
$Z_G(0,\bx)\ge\bx^E>0$, we have, for each $A\subseteq E$,
\[
\lim_{q\downarrow0}\mu_{G,q,\bx}(A)=\mu_{G,0,\bx}(A).
\]
Write
\[
p_g:=\frac{x_g}{1+x_g},\qquad X_g(A):=\ind_{\{g\in A\}},
\]
where $\ind$ denotes the indicator. Thus $\mu_{G,0,\bx}$ is independent bond
percolation with parameters $(p_g)_{g\in E}$ conditioned on connectivity.
The choice $x_g=1$ for every $g$ gives the uniform measure on connected
spanning subgraphs.

For distinct edges $e,f$, put
\begin{equation}
R_{G,q,\bx}(e,f):=
\frac{\mu_{G,q,\bx}(X_e=X_f=1)}
{\mu_{G,q,\bx}(X_e=1)\mu_{G,q,\bx}(X_f=1)}.
\label{eq:ratio}
\end{equation}
The denominator is positive, since the configuration $E$ has positive
probability. We abbreviate this ratio to $R(e,f)$ when its parameters are fixed.
Pairwise negative correlation is the assertion that $R(e,f)\le1$ for all
distinct $e,f$. For simple graphs, the uniform connected-subgraph case
is part of a conjecture of Grimmett and Winkler
\cite[Conjecture 1.1]{GrimmettWinkler2004}. The corresponding problem for
$0<q<1$ is discussed in Section 4 of the same paper.

The Lorentzian property of the homogenized multivariate Tutte polynomial,
proved by Br\"and\'en and Huh \cite[Theorem 4.10]{BrandenHuh2020}, implies
\begin{equation}
R(e,f)\le 2\left(1-\frac1m\right)
\label{eq:lorentzian-ratio}
\end{equation}
by their Proposition 2.19, with the $q=0$ case following by continuity.
Our aim is to retain information lost in this direct estimate: the
marked activities, the approach to independence as $q\uparrow1$, and,
when available, connectivity information after deleting the marked edges.
The improvement below concerns random-cluster measures; it does not
improve the universal Rayleigh constant for all Lorentzian polynomials.
To explain the coefficient comparison used below,
fix $q>0$ and the activities outside $\{e,f\}$, and write
\[
Z_G(q,\bx)=a+bx_e+cx_f+dx_ex_f,
\qquad \rho:=\frac{ad}{bc}.
\]
The coefficients $a,b,c,d$ are positive and independent of $x_e,x_f$;
their sums are given in \eqref{eq:coefficient-sums}. Direct calculation gives
\[
R(e,f)-1
=\frac{ad-bc}{(b+dx_f)(c+dx_e)}
=\frac{\rho-1}{(1+(d/b)x_f)(1+(d/c)x_e)}.
\]
In particular, $R(e,f)\le1$ if and only if $ad\le bc$. The ratio
$R(e,f)$ lies between $1$ and $\rho$: the inequality $R(e,f)\le\rho$
holds when $\rho\ge1$, whereas $R(e,f)\ge\rho$ when $\rho\le1$.
The Lorentzian estimate gives $\rho\le2(1-1/m)$ by evaluating at
$x_e=x_f=0$. We retain the dependence on the two positive marked
activities in the following refinement. Put $\gamma_m:=1-2/m$ and define
\begin{equation}
\theta_m(q):=
\begin{cases}
\displaystyle\min\left\{1-\frac2m,\frac{(1-q)^2}{4q}\right\},&0<q\le1,\\[4pt]
\displaystyle1-\frac2m,&q=0.
\end{cases}
\label{eq:theta}
\end{equation}

\begin{theorem}\label{thm:main}
	Let $G=(V,E)$ be a finite connected undirected graph with $m\ge2$ edges,
	with loops and parallel edges allowed. Let $q\in[0,1]$ and
	$\bx\in(0,\infty)^E$. For every pair of distinct edges $e,f\in E$,
	\begin{equation}
	R(e,f)\le
	1+\frac{\theta_m(q)(1-p_e)(1-p_f)}
	{(1+\theta_m(q)p_e)(1+\theta_m(q)p_f)}.
	\label{eq:main}
	\end{equation}
	In particular,
	\begin{equation}
	R(e,f)\le1+\frac{(1-p_e)(1-p_f)}{(1+p_e)(1+p_f)}<2.
	\label{eq:uniform-bound}
	\end{equation}
\end{theorem}

\begin{corollary}\label{cor:uniform}
	If $\prob$ is the uniform probability measure on connected spanning subgraphs of $G$,
	then, for distinct $e,f\in E$,
	\begin{equation}
	\begin{aligned}
	\prob(X_e=X_f=1)
	&\le\left(1+\frac{m(m-2)}{(3m-2)^2}\right)
	\prob(X_e=1)\prob(X_f=1)\\
	&\le\frac{10}{9}\prob(X_e=1)\prob(X_f=1).
	\end{aligned}
	\label{eq:ten-ninths}
	\end{equation}
\end{corollary}

\begin{proof}
	Apply Theorem~\ref{thm:main} with $q=0$ and $p_e=p_f=1/2$.
	The error term in \eqref{eq:main} becomes
	\[
	\frac{\theta_m(0)}{(2+\theta_m(0))^2}
	=\frac{m(m-2)}{(3m-2)^2}.
	\]
	The last inequality follows from
	\[
	\frac19-\frac{m(m-2)}{(3m-2)^2}
	=\frac{(3m-2)^2-9m(m-2)}{9(3m-2)^2}
	=\frac{6m+4}{9(3m-2)^2}>0\qquad(m\ge2).
	\]
\end{proof}

For $m>2$, Theorem~\ref{thm:main} strictly improves the direct Lorentzian
upper bound \eqref{eq:lorentzian-ratio} at every pair of positive marked
activities. In fact, since $\theta_m(q)\le\gamma_m$ and $0<p_e,p_f<1$,
\begin{equation}
1+\frac{\theta_m(q)(1-p_e)(1-p_f)}
{(1+\theta_m(q)p_e)(1+\theta_m(q)p_f)}
<1+\gamma_m=2\left(1-\frac1m\right).
\label{eq:strict-improvement}
\end{equation}
For $m=2$, both bounds are one. In particular, in the uniform
connected-subgraph case our factor is less than $10/9$, whereas the
direct Lorentzian factor tends to two as $m\to\infty$.
The improvement comes from combining the Lorentzian estimate $ad\leq 2(1-1/m)bc$ with  the coefficient inequalities $a\leq b$ and $a\leq c$. Using the exact expression for $R(e,f)-1$, these  inequalities allow us to retain the marked activities $x_e,x_f$ in the denominator of the error term. Since both activities are positive, the resulting upper bound is  strictly smaller  than the direct Lorentzian estimate when $m>2$. The additional $q$-dependent estimate further sharpens the bound on $\frac{ad}{bc}$ near $q=1$.

For complete graphs, Tang and Zhang
\cite[Part I, Theorem 1.4]{TangZhang2026} prove pairwise negative
correlation for the uniform connected-subgraph measure when the number
of vertices is sufficiently large. \c{C}i\c{c}eksiz and Ravichandran
\cite{CiceksizRavichandran2026} prove it on every complete graph for
common activities at least one and $0<q<1$, with the connected-subgraph
endpoint obtained by a limit. Tang and Zhang
\cite{TangZhangRegular2026} also obtain negative correlation on suitable
high-degree regular graph sequences. These results give the stronger
factor one under their respective hypotheses. Our contribution is an
explicit bound for every finite connected graph and arbitrary positive
activities, without symmetry, degree or isoperimetric assumptions.
We do not prove pairwise negative correlation on general graphs or claim
that $10/9$ is optimal.

We next refine \eqref{eq:main}. Fix $0<q\le1$ and distinct
$e,f\in E$, and let $F:=E\setminus\{e,f\}$. On subsets of $F$, define
\begin{equation}
\nu(\eta):=
\frac{q^{k(\eta)-1}\bx^\eta}
{\displaystyle\sum_{\xi\subseteq F}q^{k(\xi)-1}\bx^\xi}.
\label{eq:deleted-law}
\end{equation}
All component counts in this formula use the vertex set $V$.
This measure is defined even when $(V,F)$ is disconnected. It is the
conditional law under $\mu_{G,q,\bx}$ given $X_e=X_f=0$, not in general
the marginal law on $F$.

For $\eta\subseteq F$ and $a\in V$, let $C_\eta(a)$ be the connected
component containing $a$ in $(V,\eta)$. Write $e=\{u,v\}$ and
$f=\{w,z\}$, and define the event
\begin{equation}
J:=\left\{\eta\subseteq F:
\begin{array}{l}
C_\eta(u)\ne C_\eta(v),\quad C_\eta(w)\ne C_\eta(z),\\
\{C_\eta(u),C_\eta(v)\}=\{C_\eta(w),C_\eta(z)\}
\end{array}\right\}.
\label{eq:event-j}
\end{equation}
Thus, on $J$, each of the two edges joins the same two distinct components.
Neither $\nu$ nor $J$ depends on $x_e$ or $x_f$.
For $r\in\mathbb R$, let $[r]_+:=\max\{r,0\}$, and put
\begin{equation}
t_{e,f}(q):=
\begin{cases}
\displaystyle\left[\min\left\{1-\frac2m,
\frac{(1-q)^2}{4q}-(1-q)\nu(J)\right\}\right]_+,
&0<q\le1,\\[5pt]
\displaystyle1-\frac2m,&q=0.
\end{cases}
\label{eq:corrected-t}
\end{equation}
The dependence on $G$ and the unmarked activities is suppressed in this
notation. The value at zero is the continuous extension: since
$0\le\nu(J)\le1$,
\begin{equation}
\frac{(1-q)^2}{4q}-(1-q)\nu(J)
\ge\frac{(1-q)^2}{4q}-(1-q)\longrightarrow+\infty.
\label{eq:small-q-saturation}
\end{equation}
Consequently $t_{e,f}(q)=\theta_m(q)=\gamma_m$ for all sufficiently
small positive $q$, uniformly in the unmarked activities. No definition
of $\nu$ at zero is needed, even if deleting $e,f$ disconnects $G$.

\begin{theorem}\label{thm:correction}
	Let $G=(V,E)$ be a finite connected undirected graph with $m\ge2$ edges,
	with loops and parallel edges allowed. Let $q\in[0,1]$ and
	$\bx\in(0,\infty)^E$. For every pair of distinct edges $e,f\in E$,
	\begin{equation}
	R(e,f)\le1+
	\frac{t_{e,f}(q)(1-p_e)(1-p_f)}
	{(1+t_{e,f}(q)p_e)(1+t_{e,f}(q)p_f)}.
	\label{eq:correction-bound}
	\end{equation}
\end{theorem}

The deletion calculation in the proof of Proposition 3.2 of
\c{C}i\c{c}eksiz and Ravichandran \cite{CiceksizRavichandran2026}
relates two-edge correlations to $\nu(J)$ and a covariance under $\nu$.
We use the identity proved there to retain a nonnegative term that can
be discarded in the proof of Theorem~\ref{thm:main}. A proof of the
identity, with the normalization \eqref{eq:measure}, is included in
Section~\ref{sec:deletion}.

Theorem~\ref{thm:correction} refines Theorem~\ref{thm:main} by using
the deleted-graph probability $\nu(J)$. For $0<q<1$, it gives
pairwise negative correlation for the specified pair whenever
\begin{equation}
\nu(J)\ge\frac{1-q}{4q}.
\label{eq:sufficient}
\end{equation}
This is a sufficient condition, not an assertion that the condition
holds on every graph. The two bounds coincide at $q=0$, and both are
one at $q=1$.

Our third bound uses conditional probabilities instead of the
Lorentzian coefficient estimate. Like Theorem~\ref{thm:main}, it is
explicit in $q,p_e,p_f$ and requires no deleted-graph computation.

\begin{theorem}[Conditional-probability bound]\label{thm:conditional}
	Let $G=(V,E)$ be a finite connected undirected graph with $m\ge2$ edges,
	with loops and parallel edges allowed. Let $q\in[0,1]$ and
	$\bx\in(0,\infty)^E$. For distinct $e,f\in E$, put
	$s_g:=q+(1-q)p_g$ for $g\in\{e,f\}$. Then
	\begin{equation}
	R(e,f)\le1+
	\frac{(1-q)^2(1-p_e)(1-p_f)}{(\sqrt{s_e}+\sqrt{s_f})^2}.
	\label{eq:conditional-bound}
	\end{equation}
\end{theorem}

\begin{remark}[Comparison and use of the three bounds]\label{rem:comparison}
	Write
	\[
	H(t):=\frac{t(1-p_e)(1-p_f)}{(1+tp_e)(1+tp_f)},\quad
	B_1:=1+H(\theta_m(q)),\quad B_3:=1+H(t_{e,f}(q)),
	\]
	and let $B_4$ denote the right-hand side of
	\eqref{eq:conditional-bound}. The subscripts correspond to
	Theorems~\ref{thm:main}, \ref{thm:correction} and~\ref{thm:conditional}.
	For every admissible parameter choice,
	\begin{equation}
	0\le t_{e,f}(q)\le\theta_m(q)\le\gamma_m,\qquad
	B_3\le B_1,
	\label{eq:hierarchy}
	\end{equation}
	because $H$ is strictly increasing on $[0,1]$; see
	\eqref{eq:t-monotonicity}. Thus Theorem~\ref{thm:correction}
	implies Theorem~\ref{thm:main} as a numerical upper bound, with strict
	improvement exactly when $t_{e,f}(q)<\theta_m(q)$.
	Theorem~\ref{thm:main} is retained because it depends only on
	$m,q,p_e,p_f$ and can be evaluated immediately. In contrast,
	$\nu(J)$ depends on the graph and all unmarked activities; its exact
	definition involves sums over subsets of $E\setminus\{e,f\}$ and may
	be difficult to evaluate. When this information is unavailable, one
	can use $\nu(J)\ge0$ to recover Theorem~\ref{thm:main}.
	
	Neither $B_4$ nor $B_1$ is universally smaller. Nor is there a universal
	ordering between $B_4$ and $B_3$. Table~\ref{tab:comparison} gives
	examples on a single triangle, with $x_e=x_f=1$. In its first row,
	$B_1=B_3<B_4$; in its second row, $B_3<B_4<B_1$; in its third row,
	$B_4<B_3<B_1$. For the positive-$q$ rows,
	$\nu(J)=x_g/(q+x_g)$; the calculation is given in
	Example~\ref{ex:triangle}. These are comparisons of the upper bounds,
	not claims about equality with the actual correlation ratio.
	Combining the three theorem bounds therefore gives
	\begin{equation}
	R(e,f)\le\min\{B_3,B_4\}\qquad(0\le q\le1).
	\label{eq:combined-bound}
	\end{equation}
	Without deleted-graph information, $\min\{B_1,B_4\}$ is an entirely
	explicit alternative. At $q=0$, $B_3=B_1$ and
	$B_4=1+(1-p_e)(1-p_f)/(\sqrt{p_e}+\sqrt{p_f})^2$; at $q=1$, all
	three bounds are one.
\end{remark}

\begin{table}[htbp]
	\centering
	\renewcommand{\arraystretch}{1.35}
	\begin{tabular}{@{}ccccc@{}}
		\hline
		$q$ & $x_g$ & $B_1$ & $B_3$ & $B_4$\\
		\hline
		$0$ & $1$ & $52/49$ & $52/49$ & $9/8$\\
		$1/2$ & $1$ & $297/289$ & $1$ & $49/48$\\
		$1/2$ & $1/70$ & $297/289$ & $370/361$ & $49/48$\\
		\hline
	\end{tabular}
	\caption{Comparison on a triangle with edges $e,f,g$ and $x_e=x_f=1$.}
	\label{tab:comparison}
\end{table}

The improvement over the direct Lorentzian bound is most useful when
the marked activities stay positive, as in the connected-subgraph
limit. In the forest and spanning-tree scalings they tend to zero,
and Theorems~\ref{thm:main} and~\ref{thm:correction} recover the
Lorentzian constant rather than improving it. For spanning trees the classical negative-correlation
bound is stronger. Section~\ref{sec:models} gives these limits in detail.
Section~\ref{sec:elementary} contains the elementary inequalities,
Section~\ref{sec:deletion} proves Theorems~\ref{thm:main}
and~\ref{thm:correction}, and Section~\ref{sec:conditional} proves
Theorem~\ref{thm:conditional}.

\section{Two elementary inequalities}
\label{sec:elementary}

\begin{lemma}\label{lem:coefficient}
	Let $a,b,c,d,x,y>0$ and $0\le t\le1$. If $a\le b$, $a\le c$ and
	$ad\le(1+t)bc$, then
	\begin{equation}
	\frac{d(a+bx+cy+dxy)}{(b+dy)(c+dx)}
	\le1+\frac{t}{(1+(1+t)x)(1+(1+t)y)}.
	\label{eq:coefficient-transfer}
	\end{equation}
\end{lemma}

\begin{proof}
	Expanding the two products gives
	\[
	d(a+bx+cy+dxy)-(b+dy)(c+dx)=ad-bc.
	\]
	If $ad\le bc$, the left-hand side of \eqref{eq:coefficient-transfer}
	is at most one. Otherwise, write $\rho:=ad/(bc)$. Then
	$1<\rho\le1+t\le2$, and
	\[
	\frac dc=\rho\frac ba\ge\rho,
	\qquad \frac db=\rho\frac ca\ge\rho.
	\]
	Consequently,
	\begin{align*}
	\frac{d(a+bx+cy+dxy)}{(b+dy)(c+dx)}-1
	&=\frac{\rho-1}{(1+(d/b)y)(1+(d/c)x)}\le\frac{\rho-1}{(1+\rho x)(1+\rho y)}.
	\end{align*}
	For $1\le r\le2$,
	\begin{align*}
	\frac{d}{dr}\frac{r-1}{(1+rx)(1+ry)}
	&=\frac{(1+rx)(1+ry)-(r-1)(x+y+2rxy)}
	{(1+rx)^2(1+ry)^2}\\
	&=\frac{(1+x)(1+y)-xy(r-1)^2}{(1+rx)^2(1+ry)^2}>0,
	\end{align*}
	since the last numerator is at least $1+x+y$. Replacing $\rho$ by
	$1+t$ proves the result.
\end{proof}

For later use, recall that $p_e,p_f\in(0,1)$. Then
\begin{equation}
\frac{d}{dt}\frac{t(1-p_e)(1-p_f)}{(1+tp_e)(1+tp_f)}
=\frac{(1-p_e)(1-p_f)(1-t^2p_ep_f)}{(1+tp_e)^2(1+tp_f)^2}>0
\quad(0\le t\le1).
\label{eq:t-monotonicity}
\end{equation}

\begin{lemma}\label{lem:bounded}
	Let $U,V$ be random variables on a common probability space. If
	$1\le U\le A$ and $1\le V\le B$ almost surely, where $A,B\ge1$, then
	\begin{equation}
	\frac{\expect[UV]}{\expect U\expect V}
	\le1+\frac{(A-1)(B-1)}{(\sqrt A+\sqrt B)^2}.
	\label{eq:bounded-variables}
	\end{equation}
\end{lemma}

\begin{proof}
	If $A=1$ or $B=1$, one variable is constant and both sides are one.
	Otherwise, set
	\[
	r:=\frac{\expect U-1}{A-1},\qquad
	s:=\frac{\expect V-1}{B-1}.
	\]
	Both $r$ and $s$ belong to $[0,1]$. The pointwise inequalities
	\[
	(U-1)(V-1)\le(B-1)(U-1),\qquad
	(U-1)(V-1)\le(A-1)(V-1)
	\]
	give
	\[
	\expect[(U-1)(V-1)]\le(A-1)(B-1)\min\{r,s\}.
	\]
	Subtracting $(\expect U-1)(\expect V-1)$ and dividing by
	$\expect U\expect V$ yields
	\begin{equation}
	\frac{\expect[UV]}{\expect U\expect V}-1
	\le\frac{(A-1)(B-1)(\min\{r,s\}-rs)}
	{(1+(A-1)r)(1+(B-1)s)}.
	\label{eq:normalized-variables}
	\end{equation}
	Interchanging the variables if necessary, assume $r\le s$. Since
	\[
	\frac{d}{ds}\frac{1-s}{1+(B-1)s}
	=-\frac{B}{(1+(B-1)s)^2}<0,
	\]
	the right-hand side of \eqref{eq:normalized-variables} is at most its
	value at $s=r$. Thus
	\[
	\frac{\expect[UV]}{\expect U\expect V}-1
	\le\frac{(A-1)(B-1)r(1-r)}
	{(1+(A-1)r)(1+(B-1)r)}.
	\]
	Finally,
	\begin{align*}
	&(1+(A-1)r)(1+(B-1)r)
	-(\sqrt A+\sqrt B)^2r(1-r)\\
	&\hspace{25mm}=\bigl(1-(1+\sqrt{AB})r\bigr)^2\ge0.
	\end{align*}
	Substitution proves \eqref{eq:bounded-variables}.
\end{proof}

The bound in Lemma~\ref{lem:bounded} is attained among bounded random
variables: take $(U,V)=(A,B)$ with probability $(1+\sqrt{AB})^{-1}$
and $(U,V)=(1,1)$ otherwise. This does not assert that these extremizers
arise from a random-cluster measure.

\section{Deletion of two edges}
\label{sec:deletion}

Fix $0<q\le1$ and distinct $e,f\in E$. Set $x:=x_e$, $y:=x_f$ and
$F:=E\setminus\{e,f\}$, with the other activities held fixed. Then
\begin{equation}
Z_G(q,\bx)=a+bx+cy+dxy,
\label{eq:four-coefficients}
\end{equation}
where
\begin{equation}
\begin{aligned}
a&=\sum_{\eta\subseteq F}q^{k(\eta)-1}\bx^\eta,&
b&=\sum_{\eta\subseteq F}q^{k(\eta\cup\{e\})-1}\bx^\eta,\\
c&=\sum_{\eta\subseteq F}q^{k(\eta\cup\{f\})-1}\bx^\eta,&
d&=\sum_{\eta\subseteq F}q^{k(\eta\cup\{e,f\})-1}\bx^\eta.
\end{aligned}
\label{eq:coefficient-sums}
\end{equation}
These coefficients are positive. Adding an edge cannot increase the
number of components, so comparison of summands gives
\begin{equation}
a\le b,\qquad a\le c,\qquad b\le d,\qquad c\le d.
\label{eq:coefficient-order}
\end{equation}

\subsection{The Lorentzian bound}

\begin{lemma}\label{lem:lorentzian}
	The coefficients in \eqref{eq:coefficient-sums} satisfy
	\begin{equation}
	ad\le2\left(1-\frac1m\right)bc.
	\label{eq:lorentzian-coefficients}
	\end{equation}
\end{lemma}

\begin{proof}
	Introduce variables $z_0$ and $\bz=(z_g)_{g\in E}$, and set
	\[
	\mathcal Z(z_0,\bz):=
	\sum_{A\subseteq E}q^{k(A)-1}z_0^{m-|A|}\prod_{g\in A}z_g.
	\]
	The rank of $A$ in the graphic matroid of $G$ is
	$\rank(A)=|V|-k(A)$. Hence
	\[
	\mathcal Z(z_0,\bz)
	=q^{|V|-1}\sum_{A\subseteq E}q^{-\rank(A)}
	z_0^{m-|A|}\prod_{g\in A}z_g.
	\]
	The sum is the homogeneous multivariate Tutte polynomial of degree $m$
	used in \cite[Section 4.3]{BrandenHuh2020}. By their Theorem 4.10 it is
	Lorentzian for $0<q\le1$, and multiplication by the positive scalar
	$q^{|V|-1}$ preserves this property. Proposition 2.19 of the same paper
	therefore gives
	\begin{equation}
	\mathcal Z\,\partial_e\partial_f\mathcal Z
	\le2\left(1-\frac1m\right)
	(\partial_e\mathcal Z)(\partial_f\mathcal Z),
	\label{eq:lorentzian-derivatives}
	\end{equation}
	for all $(z_0,\bz)\in[0,\infty)^{m+1}$, where $\partial_g$ denotes
	differentiation with respect to $z_g$. This is the nonnegative-orthant
	domain in the definition of the Rayleigh property used in
	\cite[Section 2.3]{BrandenHuh2020}. Alternatively, the inequality on the
	positive orthant extends to the nonnegative orthant because both sides
	are polynomials and hence continuous.
	Evaluate at $z_0=1$, $z_e=z_f=0$ and $z_g=x_g$ for $g\in F$.
	At this point, the values of
	$\mathcal Z,\partial_e\mathcal Z,\partial_f\mathcal Z,
	\partial_e\partial_f\mathcal Z$ are $a,b,c,d$, respectively.
	Substitution proves \eqref{eq:lorentzian-coefficients}.
\end{proof}

\subsection{The deletion identity}

Recall from \eqref{eq:event-j} that $J$ is the event that, in $(V,\eta)$,
the two marked edges join the same two distinct connected components.
For $\eta\subseteq F$, define
\begin{equation}
I_e(\eta):=k(\eta)-k(\eta\cup\{e\}),\quad
I_f(\eta):=k(\eta)-k(\eta\cup\{f\}),\quad
U:=q^{-I_e},\quad V:=q^{-I_f}.
\label{eq:indicators}
\end{equation}
The variables $I_e,I_f$ take values in $\{0,1\}$. With $\nu$ as in
\eqref{eq:deleted-law}, we have
\begin{equation}
1\le U,V\le q^{-1},\qquad
\frac ba=\expect_\nu U,\qquad \frac ca=\expect_\nu V.
\label{eq:one-edge-ratios}
\end{equation}

\begin{lemma}\label{lem:deletion}
	For $0<q\le1$, the coefficients in \eqref{eq:coefficient-sums} satisfy
	\begin{equation}
	\frac da=\expect_\nu[UV]-\frac{1-q}{q^2}\nu(J).
	\label{eq:exact-deletion}
	\end{equation}
	Consequently,
	\begin{equation}
	\frac{ad}{bc}-1
	=\frac{(1-q)^2\Cov_\nu(I_e,I_f)-(1-q)\nu(J)}
	{q^2\expect_\nu U\expect_\nu V},
	\label{eq:exact-covariance}
	\end{equation}
	where $\Cov_\nu(Y_1,Y_2):=\expect_\nu[Y_1Y_2]
	-\expect_\nu Y_1\expect_\nu Y_2$.
\end{lemma}

\begin{proof}
	Contract each connected component of $(V,\eta)$ to a vertex, retaining
	the labels of edges outside $\eta$, and denote the resulting multigraph
	by $G/\eta$. An added edge reduces the component count exactly when
	its image is not a loop in $G/\eta$. If $I_e+I_f\le1$, the total
	reduction on adding both edges is $I_e+I_f$. If $I_e=I_f=1$, the images
	of both edges in $G/\eta$ are nonloops.
	They reduce the component count by one when they are parallel in $G/\eta$, which
	is precisely $J$, and by two otherwise. In the latter case they either
	share one vertex or have disjoint endpoint sets. Thus, in all cases,
	\begin{equation}
	k(\eta)-k(\eta\cup\{e,f\})
	=I_e(\eta)+I_f(\eta)-\ind_J(\eta).
	\label{eq:pointwise-count}
	\end{equation}
	On $J$ we have $UV=q^{-2}$. It follows that
	\begin{align*}
	q^{k(\eta\cup\{e,f\})-k(\eta)}
	&=UVq^{\ind_J}
	=UV\bigl(1-(1-q)\ind_J\bigr)=UV-\frac{1-q}{q^2}\ind_J.
	\end{align*}
	Taking expectation under $\nu$ gives \eqref{eq:exact-deletion}.
	Dividing by $(b/a)(c/a)$ gives
	\begin{equation}
	\frac{ad}{bc}
	=\frac{\expect_\nu[UV]}{\expect_\nu U\expect_\nu V}
	-\frac{(1-q)\nu(J)}{q^2\expect_\nu U\expect_\nu V}.
	\label{eq:normalized-deletion}
	\end{equation}
	Finally,
	\[
	U=1+\frac{1-q}{q}I_e,\qquad
	V=1+\frac{1-q}{q}I_f,
	\]
	so that $\Cov_\nu(U,V)=(1-q)^2q^{-2}\Cov_\nu(I_e,I_f)$.
	Subtracting one in \eqref{eq:normalized-deletion} proves
	\eqref{eq:exact-covariance}.
\end{proof}

The events that the endpoints of $e$ and $f$ are already connected
have indicators $1-I_e$ and $1-I_f$, whose covariance equals
$\Cov_\nu(I_e,I_f)$. Thus \eqref{eq:exact-covariance} is also the
deleted-graph formula in the proof of
\cite[Proposition 3.2]{CiceksizRavichandran2026}.

\begin{lemma}\label{lem:coefficient-bounds}
	For $0<q\le1$, the coefficients in \eqref{eq:coefficient-sums} satisfy
	\begin{equation}
	\frac{ad}{bc}\le1+t_{e,f}(q)\le1+\theta_m(q).
	\label{eq:new-coefficient-bound}
	\end{equation}
	In particular,
	\begin{equation}
	\frac{ad}{bc}\le1+\theta_m(q).
	\label{eq:old-coefficient-bound}
	\end{equation}
\end{lemma}

\begin{proof}
	Lemma~\ref{lem:bounded}, applied with $A=B=q^{-1}$, gives
	\[
	\frac{\expect_\nu[UV]}{\expect_\nu U\expect_\nu V}
	\le1+\frac{(q^{-1}-1)^2}{4q^{-1}}
	=1+\frac{(1-q)^2}{4q}.
	\]
	Substitution in \eqref{eq:normalized-deletion} gives
	\begin{equation}
	\frac{ad}{bc}-1
	\le\frac{(1-q)^2}{4q}
	-\frac{(1-q)\nu(J)}{q^2\expect_\nu U\expect_\nu V}.
	\label{eq:stronger-correction}
	\end{equation}
	By \eqref{eq:one-edge-ratios},
	$0<q^2\expect_\nu U\expect_\nu V\le1$. Therefore
	\begin{equation}
	\frac{ad}{bc}-1
	\le\frac{(1-q)^2}{4q}-(1-q)\nu(J).
	\label{eq:simpler-correction}
	\end{equation}
	Combining this with the Lorentzian bound in
	Lemma~\ref{lem:lorentzian}, and writing $C(q):=(1-q)^2/(4q)$,
	we obtain the two successive estimates
	\begin{equation}
	\begin{aligned}
	\frac{ad}{bc}-1
	&\le\min\{\gamma_m,C(q)-(1-q)\nu(J)\}\\
	&\le\bigl[\min\{\gamma_m,C(q)-(1-q)\nu(J)\}\bigr]_+
	=t_{e,f}(q)\\
	&\le\min\{\gamma_m,C(q)\}=\theta_m(q).
	\end{aligned}
	\label{eq:coefficient-hierarchy}
	\end{equation}
	The last inequality uses $\nu(J)\ge0$ and the fact that both
	$\gamma_m$ and $C(q)$ are nonnegative. This proves both assertions.
	
	This calculation identifies the loss in passing from
	Theorem~\ref{thm:correction} to Theorem~\ref{thm:main}:
	the nonnegative correction $(1-q)\nu(J)$ is discarded before the
	coefficient estimate is transferred to $R(e,f)$. The transfer function
	is increasing by \eqref{eq:t-monotonicity}, so replacing $t_{e,f}(q)$
	by $\theta_m(q)$ also weakens the final probability bound.
	The clipping at zero ensures that
	$0\le t_{e,f}(q)\le\theta_m(q)\le\gamma_m<1$ and records the
	conclusion $R(e,f)\le1$ when the corrected coefficient estimate is
	nonpositive. This is why the explicit theorem follows from the refined
	one, and also why it remains useful when $\nu(J)$ cannot be computed.
\end{proof}

\subsection{Proofs of the main results}

\begin{proof}[Proof of Theorems~\ref{thm:main} and~\ref{thm:correction}]
	For $q>0$, write $Z:=Z_G(q,\bx)$ and $\mu:=\mu_{G,q,\bx}$.
	Let $a,b,c,d$ be as in \eqref{eq:coefficient-sums}, and put
	$x:=x_e$, $y:=x_f$. By \eqref{eq:four-coefficients},
	\[
	Z=a+bx+cy+dxy.
	\]
	Consequently,
	\[
	\mu(X_e=1)=\frac{x(b+dy)}Z,\qquad
	\mu(X_f=1)=\frac{y(c+dx)}Z,\qquad
	\mu(X_e=X_f=1)=\frac{dxy}Z.
	\]
	Hence
	\begin{equation}
	R(e,f)=\frac{dZ}{(b+dy)(c+dx)}
	=1+\frac{ad-bc}{(b+dy)(c+dx)}.
	\label{eq:ratio-coefficients}
	\end{equation}
	The same expansion gives
	\[
	\Cov_\mu(X_e,X_f)=\frac{xy(ad-bc)}{Z^2}.
	\]
	Thus pairwise negative correlation for this pair is equivalent to
	$ad\le bc$. Evaluating the sums in \eqref{eq:coefficient-sums} at $q=0$
	gives the same equivalence there, since the formulas
	for $R(e,f)$ and the covariance still have positive denominators; no
	division by $a$, $b$ or $c$ is needed for this assertion.
	
	For $q>0$, set $t_1:=\theta_m(q)$ and $t_2:=t_{e,f}(q)$.
	By \eqref{eq:coefficient-order} and Lemma~\ref{lem:coefficient-bounds},
	the hypotheses of Lemma~\ref{lem:coefficient} hold with either
	$t=t_1$ or $t=t_2$. Applying that lemma gives
	\[
	R(e,f)\le1+\frac{t}{(1+(1+t)x)(1+(1+t)y)}.
	\]
	Since $x=p_e/(1-p_e)$ and $y=p_f/(1-p_f)$, we have
	\[
	1+(1+t)x=\frac{1+tp_e}{1-p_e},\qquad
	1+(1+t)y=\frac{1+tp_f}{1-p_f}
	\]
	and hence
	\[
	R(e,f)\le1+\frac{t(1-p_e)(1-p_f)}{(1+tp_e)(1+tp_f)}.
	\]
	Taking $t=t_2$ proves the refined bound
	\eqref{eq:correction-bound}. Since $t_2\le t_1$ and the right-hand
	side is increasing in $t$, this also proves \eqref{eq:main}.
	Equivalently, the latter follows directly by taking $t=t_1$.
	Thus the implication between the coefficient bounds in
	\eqref{eq:coefficient-hierarchy} is preserved for the probability bounds.
	Since $\theta_m(q)\le1$, \eqref{eq:t-monotonicity} gives
	\eqref{eq:uniform-bound}. Its strict upper bound by two follows from
	$(1-p_e)(1-p_f)<(1+p_e)(1+p_f)$.
	
	It remains to prove both theorems at $q=0$. Keep $G$ and all
	activities fixed. The polynomial $Z_G(q,\bx)$ has nonnegative integer
	powers of $q$, and
	\[
	Z_G(0,\bx)=\sum_{A:\,k(A)=1}\bx^A\ge\bx^E>0.
	\]
	The numerators of the three probabilities in \eqref{eq:ratio} are also
	polynomials in $q$. Each marginal numerator at zero is at least
	$\bx^E$. Thus the probabilities and their ratio converge to those of
	$\mu_{G,0,\bx}$ as $q\downarrow0$. Moreover,
	\[
	\lim_{q\downarrow0}\theta_m(q)
	=\lim_{q\downarrow0}t_{e,f}(q)
	=1-\frac2m=\theta_m(0)=t_{e,f}(0),
	\]
	where the limit for $t_{e,f}$ follows from
	\eqref{eq:small-q-saturation}. Taking limits in
	\eqref{eq:main} and \eqref{eq:correction-bound} proves both assertions.
	In fact, their right-hand sides agree for all sufficiently small
	positive $q$. No quotient of the coefficients $a,b,c,d$ at zero,
	or deleted-graph law at zero, is used.
\end{proof}

\begin{remark}\label{rem:stronger}
	One can retain the denominator in \eqref{eq:stronger-correction}.
	More precisely, for $0<q\le1$, Theorem~\ref{thm:correction} remains valid if
	$t_{e,f}(q)$ is replaced by
	\[
	\left[\min\left\{1-\frac2m,\frac{(1-q)^2}{4q}
	-\frac{(1-q)\nu(J)}{q^2\expect_\nu U\expect_\nu V}
	\right\}\right]_+.
	\]
	This parameter is no larger than $t_{e,f}(q)$. Both corrections require
	information about the deleted graph $(V,E\setminus\{e,f\})$.
	We make no claim that $10/9$ is
	the best graph-independent constant.
\end{remark}

\begin{example}\label{ex:triangle}
	Let $G$ be a triangle with edges $e=\{u,v\}$, $f=\{v,w\}$ and
	$g=\{u,w\}$. After deleting $e,f$, the two configurations are
	$\varnothing$ and $\{g\}$, with weights $q^2$ and $qx_g$ under
	\eqref{eq:deleted-law}. Only the second belongs to $J$. Hence
	\[
	\nu(J)=\frac{qx_g}{q^2+qx_g}=\frac{x_g}{q+x_g}.
	\]
	If $q=1/2$ and $x_g=1$, then $\nu(J)=2/3$ and
	\[
	\theta_3(1/2)=\frac18,\qquad
	t_{e,f}(1/2)=\left[\min\left\{\frac13,
	\frac18-\frac13\right\}\right]_+=0.
	\]
	For every $x_e,x_f>0$, \eqref{eq:correction-bound} therefore gives
	$R(e,f)\le1$, whereas the right-hand side of \eqref{eq:main} is
	strictly greater than one.
\end{example}

\section{Proof of the conditional-probability bound}
\label{sec:conditional}

\begin{proof}[Proof of Theorem~\ref{thm:conditional}]
	First suppose $q>0$, and write $\mu:=\mu_{G,q,\bx}$ and
	$F:=E\setminus\{e,f\}$. Let $\omega$ be a random edge set with law
	$\mu$, and let $\bar\mu$ be the law of $\omega\cap F$.
	In contrast to \eqref{eq:deleted-law}, this measure
	is obtained by summing over the states of $e$ and $f$.
	For a fixed $\eta\subseteq F$, set
	\[
	\begin{aligned}
	a_\eta&:=q^{k(\eta)-1},&
	b_\eta&:=q^{k(\eta\cup\{e\})-1},\\
	c_\eta&:=q^{k(\eta\cup\{f\})-1},&
	d_\eta&:=q^{k(\eta\cup\{e,f\})-1},
	\end{aligned}
	\]
	and $Z_\eta:=a_\eta+b_\eta x_e+c_\eta x_f+d_\eta x_ex_f$.
	Conditional on $\omega\cap F=\eta$, the weights of the four states of the
	marked edges are $a_\eta,b_\eta x_e,c_\eta x_f,d_\eta x_ex_f$.
	Equation~\eqref{eq:pointwise-count} gives
	\[
	k(\eta)+k(\eta\cup\{e,f\})
	-k(\eta\cup\{e\})-k(\eta\cup\{f\})=\ind_J(\eta)\ge0.
	\]
	Since $q\le1$, it follows that $a_\eta d_\eta\le b_\eta c_\eta$.
	The conditional covariance is therefore
	\begin{align*}
	\Cov_\mu(X_e,X_f\mid\omega\cap F=\eta)
	&=\frac{x_ex_f}{Z_\eta^2}
	\bigl[d_\eta Z_\eta
	-(b_\eta+d_\eta x_f)(c_\eta+d_\eta x_e)\bigr]=\frac{x_ex_f(a_\eta d_\eta-b_\eta c_\eta)}{Z_\eta^2}\le0.
	\end{align*}
	Define $Y_g(\eta):=\mu(X_g=1\mid\omega\cap F=\eta)$ for $g=e,f$.
	Taking expectation of the preceding inequality yields
	\begin{equation}
	\mu(X_e=X_f=1)\le\expect_{\bar\mu}[Y_eY_f],\qquad
	\mu(X_g=1)=\expect_{\bar\mu}Y_g.
	\label{eq:conditional-mixture}
	\end{equation}
	
	For a single edge $g$ and $\xi\subseteq E\setminus\{g\}$, condition
	on $\omega\cap(E\setminus\{g\})=\xi$. If its endpoints are connected in
	$(V,\xi)$, the conditional odds of $X_g=1$ are $x_g$; otherwise
	they are $x_g/q$. Thus
	\begin{equation}
	\mu(X_g=1\mid\omega\cap(E\setminus\{g\})=\xi)=
	\begin{cases}
	\displaystyle\frac{x_g}{1+x_g}=p_g,&\text{endpoints connected},\\[5pt]
	\displaystyle\frac{x_g}{q+x_g}=\frac{p_g}{s_g},&\text{otherwise}.
	\end{cases}
	\label{eq:single-edge}
	\end{equation}
	Averaging over the remaining marked edge, conditional on $\omega\cap F$,
	shows that
	\[
	p_g\le Y_g\le\frac{p_g}{s_g}\qquad(g=e,f).
	\]
	Apply Lemma~\ref{lem:bounded} under $\bar\mu$ to
	$Y_e/p_e,Y_f/p_f$, with bounds $1/s_e,1/s_f$. By
	\eqref{eq:conditional-mixture},
	\begin{align*}
	R(e,f)
	&\le\frac{\expect_{\bar\mu}[Y_eY_f]}
	{\expect_{\bar\mu}Y_e\expect_{\bar\mu}Y_f}\le1+\frac{(s_e^{-1}-1)(s_f^{-1}-1)}
	{(s_e^{-1/2}+s_f^{-1/2})^2}=1+\frac{(1-s_e)(1-s_f)}{(\sqrt{s_e}+\sqrt{s_f})^2}.
	\end{align*}
	Since $1-s_g=(1-q)(1-p_g)$, this is
	\eqref{eq:conditional-bound} for $q>0$. At $q=0$, the probability
	ratio is continuous by the argument in Section~\ref{sec:deletion},
	and $s_g\to p_g>0$. Taking the limit proves the remaining case.
\end{proof}

For equal marked parameters $p_e=p_f=p$, the bound reads
\[
R(e,f)\le1+\frac{(1-q)^2(1-p)^2}{4(q+(1-q)p)}.
\]
For example, it gives $49/48$ at $q=p=1/2$, independently of
$m$ and the unmarked activities. The three bounds are compared in
Remark~\ref{rem:comparison} and Table~\ref{tab:comparison}.
Unlike the proof of Theorem~\ref{thm:main}, this proof uses no
Lorentzian inequality: the ingredients are conditional negative
correlation of the marked pair given all other edges, the single-edge
probability interval, and Lemma~\ref{lem:bounded}.

\begin{remark}
	If $e$ is a loop, then $k(A\cup\{e\})=k(A)$ for every
	$A\subseteq E\setminus\{e\}$, and its state is independent of all
	other edges with open probability $p_e$. If $e$ is a bridge of $G$,
	its endpoints are disconnected in every such $A$, so for $q>0$ its
	state is independent of all other edges with open probability
	$x_e/(q+x_e)$. At $q=0$ it is present in every connected spanning
	subgraph. In either case $R(e,f)=1$ for every $f\ne e$.
\end{remark}

\section{Specializations and limits}
\label{sec:models}

All random-cluster parameter limits in this section are taken on a fixed
graph $G$. Write
$N:=|V|$ and recall $\gamma_m=1-2/m$. Let $\mathcal C$, $\mathcal F$ and
$\mathcal T$ be the sets of connected spanning subgraphs, spanning
forests and spanning trees of $G$, respectively, identified with their
edge sets. Forests have no designated roots, and isolated vertices are
included. For positive weights
$\boldsymbol\lambda=(\lambda_g)_{g\in E}$, put
\begin{equation}
\mu_{\mathcal A,\boldsymbol\lambda}(A)
:=\frac{\ind_{\{A\in\mathcal A\}}\boldsymbol\lambda^A}
{\displaystyle\sum_{B\in\mathcal A}\boldsymbol\lambda^B},
\qquad \mathcal A\in\{\mathcal C,\mathcal F,\mathcal T\},
\label{eq:restricted-measures}
\end{equation}
where $\boldsymbol\lambda^A:=\prod_{g\in A}\lambda_g$.
Write $R_{\mathcal A,\boldsymbol\lambda}(e,f)$ for the ratio
\eqref{eq:ratio} with this measure in place of $\mu_{G,q,\bx}$.
For forests and trees, the marked edges are assumed not to be loops,
so that both marginal probabilities are positive. The limits below are
the weighted versions of those in
\cite[Section 4]{GrimmettWinkler2004}.

\subsection{Independent percolation}

At $q=1$,
\[
Z_G(1,\bx)=\prod_{g\in E}(1+x_g),\qquad
\mu_{G,1,\bx}(A)=\prod_{g\in A}p_g\prod_{g\notin A}(1-p_g).
\]
Thus $R(e,f)=1$, which agrees with Theorem~\ref{thm:main} since
$\theta_m(1)=t_{e,f}(1)=0$; the other two theorem bounds also equal one.

\subsection{Connected spanning subgraphs}

Keep $x_g=\lambda_g$ fixed and let $q\downarrow0$. Termwise convergence
in the finite partition sum gives
\[
Z_G(q,\boldsymbol\lambda)\longrightarrow
\sum_{A\in\mathcal C}\boldsymbol\lambda^A>0,
\qquad
\mu_{G,q,\boldsymbol\lambda}(A)\longrightarrow
\mu_{\mathcal C,\boldsymbol\lambda}(A).
\]
Using $p_g=\lambda_g/(1+\lambda_g)$ in Theorem~\ref{thm:main}, we obtain
\begin{equation}
\begin{aligned}
R_{\mathcal C,\boldsymbol\lambda}(e,f)
&\le1+\frac{\gamma_m}
{[1+(1+\gamma_m)\lambda_e][1+(1+\gamma_m)\lambda_f]} \le1+\frac{1}{(1+2\lambda_e)(1+2\lambda_f)}.
\end{aligned}
\label{eq:connected-specialization}
\end{equation}
Indeed,
\[
\frac{\gamma_m(1-p_e)(1-p_f)}
{(1+\gamma_m p_e)(1+\gamma_m p_f)}
=\frac{\gamma_m}
{[1+(1+\gamma_m)\lambda_e][1+(1+\gamma_m)\lambda_f]},
\]
and the last inequality in \eqref{eq:connected-specialization} follows
from \eqref{eq:t-monotonicity} and $0\le\gamma_m\le1$.
For a common weight $\lambda>0$, this reads
\begin{equation}
R_{\mathcal C,\lambda}(e,f)
\le1+\frac{\gamma_m}{[1+(1+\gamma_m)\lambda]^2}
\le1+\frac{1}{(1+2\lambda)^2},
\label{eq:connected-equal-weights}
\end{equation}
where a scalar subscript denotes the corresponding constant weight vector.
In particular, $\lambda=1$ gives the uniform connected-subgraph measure
and the bound in Corollary~\ref{cor:uniform}. The graph-independent
bounds at $\lambda=1/2,1,2$ are $5/4$, $10/9$ and $26/25$, respectively.

\subsection{The arboreal gas and the uniform forest}

The measure $\mu_{\mathcal F,\boldsymbol\lambda}$ is the arboreal gas
with edge weights $\boldsymbol\lambda$. With a common weight $\lambda$,
it assigns probability proportional to $\lambda^{|A|}$ to each forest;
see \cite[Section 1.1]{BauerschmidtCrawfordHelmuthSwan2021}.
For $A\subseteq E$, let
\[
\ell(A):=|A|-N+k(A)
\]
be its cycle rank. Choosing a spanning tree in each component shows that
$\ell(A)\ge0$, with equality exactly when $A\in\mathcal F$.
Set $x_g=q\lambda_g$. Then
\begin{equation}
\begin{aligned}
q^{-(N-1)}q^{k(A)-1}\prod_{g\in A}(q\lambda_g)
&=q^{|A|+k(A)-N}\boldsymbol\lambda^A
=q^{\ell(A)}\boldsymbol\lambda^A,\\
q^{-(N-1)}Z_G(q,q\boldsymbol\lambda)
&=\sum_{A\subseteq E}q^{\ell(A)}\boldsymbol\lambda^A
\longrightarrow\sum_{A\in\mathcal F}\boldsymbol\lambda^A>0.
\end{aligned}
\label{eq:forest-limit}
\end{equation}
The same normalization applies to each configuration weight, so
$\mu_{G,q,q\boldsymbol\lambda}$ converges to
$\mu_{\mathcal F,\boldsymbol\lambda}$. For either marked edge $g$,
the singleton $\{g\}$ is a forest of positive weight. Both marginal
probabilities therefore have positive limits.

Along this scaling, the percolation parameters in Theorem~\ref{thm:main}
are
\[
p_g(q)=\frac{q\lambda_g}{1+q\lambda_g}\longrightarrow0,
\qquad \theta_m(q)\longrightarrow\gamma_m.
\]
Since $R_{\mathcal F,\boldsymbol\lambda}(e,f)
=\lim_{q\downarrow0}R_{G,q,q\boldsymbol\lambda}(e,f)$, we have 
\begin{equation}
\begin{aligned}
R_{\mathcal F,\boldsymbol\lambda}(e,f)
&\le1+\lim_{q\downarrow0}
\frac{\theta_m(q)(1-p_e(q))(1-p_f(q))}
{(1+\theta_m(q)p_e(q))(1+\theta_m(q)p_f(q))}=1+\gamma_m=2\left(1-\frac1m\right).
\end{aligned}
\label{eq:forest-bound}
\end{equation}
This recovers the Lorentzian bound, rather than improving it in the
forest limit. Taking every $\lambda_g=1$ gives the uniform measure on
\emph{all} spanning forests, with no restriction on their number of
components. The same bound \eqref{eq:forest-bound} applies to this
uniform forest measure.

The bound in Theorem~\ref{thm:correction}, with the parameter as stated
there, has the same limit. Indeed, for every deleted-graph measure,
\[
\frac{(1-q)^2}{4q}-(1-q)\nu(J)
\ge\frac{(1-q)^2}{4q}-(1-q)\longrightarrow+\infty.
\]
Hence $t_{e,f}(q)=\gamma_m$ for all sufficiently small $q$.
In particular, one cannot obtain $10/9$ for the uniform forest by
substituting $p_e=p_f=1/2$: the parameters tending to zero above are
those of the random-cluster measure before taking the forest limit.

\subsection{Spanning trees}

Fix $\alpha>1$. Let $0<t<1$, take $q=t^\alpha$ and
$x_g=t\lambda_g$, and let $t\downarrow0$.
For each $A\subseteq E$,
\begin{equation}
\begin{aligned}
t^{-(N-1)}q^{k(A)-1}\prod_{g\in A}(t\lambda_g)
&=t^{\alpha(k(A)-1)+|A|-(N-1)}\boldsymbol\lambda^A\\
&=t^{(\alpha-1)(k(A)-1)+\ell(A)}\boldsymbol\lambda^A.
\end{aligned}
\label{eq:tree-limit}
\end{equation}
Since $\alpha-1>0$, $k(A)-1\ge0$ and $\ell(A)\ge0$, the exponent
is nonnegative, and it is zero exactly when $k(A)=1$
and $\ell(A)=0$, that is, when $A\in\mathcal T$.
Since $G$ is connected, $\sum_{A\in\mathcal T}\boldsymbol\lambda^A>0$.
Dividing the normalized weights in \eqref{eq:tree-limit} by their sum
therefore proves convergence to $\mu_{\mathcal T,\boldsymbol\lambda}$.
Each nonloop edge extends to a spanning tree, so the marked marginals
again have positive limits. The choice $\lambda_g=1$ gives the uniform
spanning tree (UST).

Theorem~\ref{thm:main} gives only
$R_{\mathcal T,\boldsymbol\lambda}(e,f)\le2(1-1/m)$ along this scaling,
since $p_g(t)=t\lambda_g/(1+t\lambda_g)\to0$ and
$\theta_m(t^\alpha)\to\gamma_m$. This limit is the same for every
fixed $\alpha>1$. Theorem~\ref{thm:correction}, with $t_{e,f}$ as
stated, gives the same constant: since $0\le\nu(J)\le1$,
\[
\frac{(1-t^\alpha)^2}{4t^\alpha}-(1-t^\alpha)\nu(J)
\ge\frac{(1-t^\alpha)^2}{4t^\alpha}-(1-t^\alpha)
\longrightarrow+\infty,
\]
so $t_{e,f}(t^\alpha)=\gamma_m$ for all sufficiently small $t$.
The bound in Theorem~\ref{thm:conditional} does not improve this
limit either. Set $s_g(t):=t^\alpha+(1-t^\alpha)p_g(t)$ for $g=e,f$.
Then $t^\alpha/t\to0$ and
\[
\frac{p_g(t)}t=\frac{\lambda_g}{1+t\lambda_g}\longrightarrow\lambda_g,
\qquad
\frac{s_g(t)}t=t^{\alpha-1}+(1-t^\alpha)\frac{p_g(t)}t
\longrightarrow\lambda_g.
\]
Consequently,
\[
\frac{(\sqrt{s_e(t)}+\sqrt{s_f(t)})^2}{t}
\longrightarrow(\sqrt{\lambda_e}+\sqrt{\lambda_f})^2,
\qquad
(1-t^\alpha)^2(1-p_e(t))(1-p_f(t))\longrightarrow1,
\]
and the error term on the right-hand side of
\eqref{eq:conditional-bound} tends to $+\infty$.
At $\alpha=1$ the exponent in \eqref{eq:tree-limit} is $\ell(A)$,
so the limit is instead the arboreal gas with weights
$\boldsymbol\lambda$, as in \eqref{eq:forest-limit}.

For the UST the classical bound is instead $ R_{\mathcal T,1}(e,f)\le1$; for instance see Chapter 4 of \cite{LP2016}. Thus the limiting bound supplied by Theorem~\ref{thm:main} need not be sharp even when the limiting measure is known to be negatively
correlated.


\section*{Statement of AI use}
The authors used ChatGPT (6.1 Sol) in developing and checking the proofs, as well as for improving the manuscript's writing. The authors take full responsibility for all proofs, references, and other content of this paper.


\bibliographystyle{plain}
\bibliography{references}

\end{document}